\documentclass{article}

\newcommand{\KP}{{\sf KP}}
\newcommand{\sfqk}{{\sf qk}}

\usepackage{proof}
\usepackage{latexsym}
\usepackage{amsfonts}
\usepackage{amssymb}
\usepackage{cite}

\newcommand{\blem}{\begin{lemma}}
\newcommand{\elem}{\end{lemma}}
\newcommand{\bth}{\begin{theorem}}
\newcommand{\ethm}{\end{theorem}}
\newcommand{\benu}{\begin{enumerate}}
\newcommand{\eenu}{\end{enumerate}}
\newcommand{\bdes}{\begin{description}}
\newcommand{\edes}{\end{description}}
\newcommand{\bdf}{\begin{definition}}
\newcommand{\edf}{\end{definition}}
\newcommand{\bcor}{\begin{cor}}
\newcommand{\ecor}{\end{cor}}
\newcommand{\bprp}{\begin{proposition}}
\newcommand{\eprp}{\end{proposition}}
\newcommand{\bmlem}{\begin{mlemma}}
\newcommand{\emlem}{\end{mlemma}}
\newcommand{\bclm}{\begin{claim}}
\newcommand{\eclm}{\end{claim}}
\newcommand{\bprf}{{\bf Proof}.\hspace{2mm}}

\newcommand{\eprf}{\hspace*{\fill} $\Box$}

\newcommand{\beqn}{\begin{equation}}
\newcommand{\eeqn}{\end{equation}}
\newcommand{\beqnarr}{\begin{eqnarray}}
\newcommand{\eeqnarr}{\end{eqnarray}}
\newcommand{\beqnarrs}{\begin{eqnarray*}}
\newcommand{\eeqnarrs}{\end{eqnarray*}}

\newcommand{\spand}{\,\&\,}

\newcommand{\restrict}{\!\upharpoonright\!}

\newtheorem{theorem}{Theorem}[section]
\newtheorem{definition}[theorem]{Definition}
\newtheorem{proposition}[theorem]{Proposition}
\newtheorem{lemma}[theorem]{Lemma}
\newtheorem{cor}[theorem]{Corollary}
\newtheorem{mlemma}[theorem]{Main Lemma}
\newtheorem{claim}[theorem]{Claim}

\newcommand{\alp}{\alpha}

\newcommand{\veps}{\varepsilon}
\newcommand{\del}{\delta}
\newcommand{\Del}{\Delta}
\newcommand{\ome}{\omega}

\newcommand{\bet}{\beta}
\newcommand{\gam}{\gamma}
\newcommand{\Gam}{\Gamma}
\newcommand{\kap}{\kappa}
\newcommand{\sig}{\sigma}
\newcommand{\Sig}{\Sigma}
\newcommand{\tht}{\theta}
\newcommand{\Tht}{\Theta}
\newcommand{\lam}{\lambda}
\newcommand{\Lam}{\Lambda}
\newcommand{\vphi}{\varphi}

\newcommand{\fal}{\forall}
\newcommand{\exi}{\exists}

\newcommand{\Rarw }{\Rightarrow}

\newcommand{\lrarw}{\leftrightarrow}
\newcommand{\Lrarw}{\Leftrightarrow}

\newcommand{\darw}{\downarrow}

\newcommand{\calh}{{\cal H}}

\newcommand{\calL}{{\cal L}}

\newcommand{\calP}{{\cal P}}

\newcommand{\la}{\langle}
\newcommand{\ra}{\rangle}

\newcommand{\rk}{\mbox{{\rm rk}}}

\newcommand{\sfk}{{\sf k}}

\title{
$\ome_{1}$ under $\Pi_{1}$-Collection
}
\author{Toshiyasu Arai
\\
Graduate School of Science,
Chiba University
\\
1-33, Yayoi-cho, Inage-ku,
Chiba, 263-8522, JAPAN
\\
tosarai@faculty.chiba-u.jp
}
\date{}
\begin{document}
\maketitle

\begin{abstract}
We describe a proof-theoretic bound on $\Sig_{2}$-definable countable ordinals
in Kripke-Platek set theory with $\Pi_{1}$-Collection and the existence of $\ome_{1}$.
\end{abstract}

\section{Introduction}

Let $(\ome_{1})$ denote an axiom stating that `there exists an uncountable regular ordinal', and
$T_{1}:=\KP\ome+(V=L)+(\Pi_{1}\mbox{{\rm -Collection}})+(\ome_{1})$.
Let $\rho_{0}$ denote the least ordinal above $\ome_{1}$ such that $L_{\rho_{0}}\models (\Pi_{1}\mbox{{\rm -Collection}})$.
In this note a collapsing function
$\Psi_{\ome_{1}}:\alp\mapsto\Psi_{\ome_{1}}(\alp)<\ome_{1}$ is defined, and it is shown that for each $n<\ome$,
$T_{1}\vdash\fal\alp<\ome_{n}(\rho_{0}+1)\exi x<\ome_{1}(x=\Psi_{\ome_{1}}(\alp))$ with a $\Sig_{2}$-formula
$x=\Psi_{\ome_{1}}(\alp)$, cf. Lemma \ref{lem:lowerbndreg}.
Conversely we show the
\bth\label{th:mainthZ}

For a sentence $\exi x\in L_{\ome_{1}}\, \vphi(x)$ with a $\Sig_{2}$-formula $\vphi(x)$, 
if
\[
T_{1}\vdash\exi x\in L_{\ome_{1}}\,\vphi(x) 
\]
then
\[
\exi n<\ome[T_{1}\vdash \exi x\in L_{\Psi_{\ome_{1}}(\ome_{n}(\rho_{0}+1))} \vphi(x)]
.\]
\end{theorem}

This paper relies on our \cite{liftupZF}.

\section{$\Sig_{1}$-Skolem hulls}\label{sect:Skolemhull}

Everything in this section is reproduced from \cite{liftupZF}.

For a model $\la M;\in\restrict (M\times M)\ra$ and $X\subset M$,
$\Sig_{1}^{M}(X)$ denotes the set of $\Sig_{1}(X)$-definable subsets of $M$,
where $\Sig_{1}(X)$-formulae may have parameters from $X$.
$\Sig_{1}^{M}(M)$ is denoted $\Sig_{1}(M)$.

An ordinal $\alp>1$ is said to be a {\it multiplicative principal number\/} iff
$\alp$ is closed under ordinal multiplication, i.e., $\exi \bet[\alp=\ome^{\ome^{\bet}}]$.
If $\alp$ is a multiplicative principal number, then $\alp$ is closed under G\"odel's pairing function $j$
and there exists a $\Del_{1}$-bijection between $\alp$ and $L_{\alp}$
for the constructible hierarchy $L_{\alp}$ up to $\alp$.
In this section 
$\sig$ is assumed to be a multiplicative principal number$>\ome$.

\bdf\label{df:crdS} 

\benu

\item 
$
cf(\kap) := \min\{\alp\leq\kap : \mbox{{\rm there is a cofinal map }} f:\alp\to \kap\}
$.

\item
$\rho(L_{\sig})$ {\rm denotes the} $\Sig_{1}$-projectum {\rm of} $L_{\sig}${\rm :}
$\rho(L_{\sig})$ {\rm is the least ordinal} $\rho$ {\rm such that}
$
\calP(\rho)\cap\Sig_{1}(L_{\sig})\not\subset L_{\sig}
$.


\item
{\rm Let} $\alp\leq\bet$ {\rm and} $f: L_{\alp}\to L_{\bet}$.
{\rm Then the map} $f$ {\rm is a} $\Sig_{1}$-elementary embedding{\rm , denoted}
$
f: L_{\alp}\prec_{\Sig_{1}} L_{\bet}
$
{\rm iff for any} $\Sig_{1}(L_{\alp})${\rm -sentence} $\vphi[\bar{a}]\, (\bar{a}\subset L_{\alp})$,
$
L_{\alp} \models \vphi[\bar{a}] \Lrarw 
L_{\bet} \models \vphi[f(\bar{a})]
$
{\rm where} $f(\bar{a})=f(a_{1}),\ldots,f(a_{k})$ {\rm for} $\bar{a}=a_{1},\ldots,a_{k}$.
{\rm An ordinal} $\gam$ {\rm such that}
$
\fal\del<\gam[f(\del)=\del]\spand f(\gam)>\gam
$
{\rm is said to be the} critical point {\rm of the} $\Sig_{1}$-
{\rm elementary embedding} $f$
{\rm if such an ordinal} $\gam$ {\rm exists.}

\item
{\rm For} $X\subset L_{\sig}$,
$\mbox{{\rm Hull}}^{\sig}_{\Sig_{1}}(X)$ {\rm denotes the set} 
{\rm (}$\Sig_{1}$-Skolem hull {\rm of} $X$ {\rm in} $L_{\sig}${\rm ) defined as follows.}
$<_{L}$ {\rm denotes a} $\Del_{1}${\rm -well ordering of the constructible universe} $L$.
{\rm Let} $\{\vphi_{i}:i\in\ome\}$ {\rm denote an enumeration of} $\Sig_{1}${\rm -formulae in the language}
$\{\in\}${\rm . Each is of the form} $\vphi_{i}\equiv\exi y\theta_{i}(x,y;u)\, (\tht\in\Del_{0})$ {\rm with fixed variables} $x,y,u${\rm . Set for} $b\in X$
\beqnarrs
r_{\Sig_{1}}^{\sig}(i,b) & \simeq & \mbox{ {\rm the }} <_{L} \mbox{{\rm -least }} c\in L_{\sig}
\mbox{ {\rm such that} } L_{\sig}\models\theta_{i}((c)_{0},(c)_{1}; b)
\\
h_{\Sig_{1}}^{\sig}(i,b) & \simeq & (r_{\Sig_{1}}^{\sig}(i,b))_{0}
\\
\mbox{{\rm Hull}}^{\sig}_{\Sig_{1}}(X) & = & rng(h_{\Sig_{1}}^{\sig})=\{h_{\Sig_{1}}^{\sig}(i,b)\in L_{\sig}:i\in\ome, b\in X\}
\eeqnarrs
{\rm Then}
$
L_{\sig}\models \exi x\exi y\, \theta_{i}(x,y;b) \to h_{\Sig_{1}}^{\sig}(i,b)\darw \spand \exi y\, \theta_{i}(h_{\Sig_{1}}^{\sig}(i,b),y;b)
$.
\eenu
\edf

\bprp\label{prp:rhHull}
Assume that $X$ is a set in $L_{\sig}$.
Then
$r_{\Sig_{1}}^{\sig}$ and $h_{\Sig_{1}}^{\sig}$ are partial $\Del_{1}(L_{\sig})$-maps such that the domain of $h_{\Sig_{1}}^{\sig}$ is a 
$\Sig_{1}(L_{\sig})$-subset of 
$\ome\times X$.
Therefore its range $\mbox{{\rm Hull}}^{\sig}_{\Sig_{1}}(X)$ is a $\Sig_{1}(L_{\sig})$-subset of $L_{\sig}$.
%
\eprp


\bprp\label{clm:crdlocal3.1} 
Let $Y=\mbox{{\rm Hull}}^{\sig}_{\Sig_{1}}(X)$.
For any $\Sig_{1}(Y)$-sentence $\vphi(\bar{a})$ with parameters $\bar{a}$ from $Y$
$
L_{\sig}\models \vphi(\bar{a})\Lrarw Y\models \vphi(\bar{a})
$.
Namely
$
Y \prec_{\Sig_{1}}L_{\sig}
$.
\eprp

\bdf\label{df:pikap}{\rm (Mostowski collapsing function} $F${\rm )}

{\rm By Proposition \ref{clm:crdlocal3.1} and the Condensation Lemma
we have an isomorphism (Mostowski collapsing function) }
\[
F:\mbox{{\rm Hull}}^{\sig}_{\Sig_{1}}(X)\lrarw  L_{\gam}
\]
{\rm for an ordinal} $\gam\leq\sig$ {\rm such that} $F\restrict Y=id\restrict Y$ 
{\rm for any transitive}
$Y\subset \mbox{{\rm Hull}}^{\sig}_{\Sig_{1}}(X)$.

{\rm Let us denote, though} $\sig\not\in dom(F)=\mbox{{\rm Hull}}^{\sig}_{\Sig_{1}}(X)$
\[
F(\sig):=\gam
.\]
{\rm Also for the above Mostowski collapsing map} $F$ {\rm let}
\[
F^{\Sig_{1}}(x;\sig,X):=F(x)
.\]
{\rm The inverse} $G:=F^{-1}$ {\rm of} $F$ {\rm is a} $\Sig_{1}${\rm -elementary embedding
from} $L_{F(\sig)}$ {\rm to} $L_{\sig}$.
\edf

\bprp\label{prp:complexcr}
Let 
$L_{\sig}\models\KP\ome+\Sig_{1}\mbox{{\rm -Collection}}$.
Then for $\kap\leq\sig$,
$\{(x,y): x<\kap\spand y=\min\{y<\kap: \mbox{{\rm Hull}}^{\sig}_{\Sig_{1}}(x\cup\{\kap\})\cap\kap\subset y\}\}$ 
is a $Bool(\Sig_{1}(L_{\sig}))$-predicate on $\kap$, and hence a set in $L_{\sig}$
if $\kap<\sig$ and $L_{\sig}\models\Sig_{1}\mbox{{\rm -Separation}}$.
\eprp


$F^{\Sig_{1}}_{x\cup\{\kap\}}(y)$ denotes the Mostowski collapse $F^{\Sig_{1}}(y;\sig,x\cup\{\kap\})$.

\bth\label{th:cofinalitylocal}
Let $\sig$ be an ordinal such that $L_{\sig}\models \mbox{{\rm KP}}\ome+\Sig_{1}\mbox{{\rm -Separation}}$,
 and $\ome\leq\alp<\kap<\sig$ with
$\alp$ a multiplicative principal number and $\kap$ a limit ordinal.
Then the following conditions are mutually equivalent:
\benu

\item 
$L_{\sig}\models {}^{\alp}\kap\subset L_{\kap}$.

\item 
$L_{\sig}\models \alp<cf(\kap)$.

\item
There exists an ordinal $x$ such that $\alp<x<\kap$,
$ \mbox{{\rm Hull}}^{\sig}_{\Sig_{1}}(x\cup\{\kap\})\cap \kap\subset x$ and $F^{\Sig_{1}}_{x\cup\{\kap\}}(\sig)<\kap$.

\item 
For the Mostowski collapse
$F^{\Sig_{1}}_{x\cup\{\kap\}}(y)$, there exists an ordinal $x$ such that 
$\alp<x=F^{\Sig_{1}}_{x\cup\{\kap\}}(\kap)<F^{\Sig_{1}}_{x\cup\{\kap\}}(\sig)<\kap$,
and for any $\Sig_{1}$-formula $\vphi$ and any $a\in L_{x}$,
$L_{\sig}\models\vphi[\kap,a]\to  L_{F^{\Sig_{1}}_{x\cup\{\kap\}}(\sig)} \models\vphi[x,a]$ holds.

\eenu

\end{theorem}

\bdf\label{df:KPome1}
$T_{1}:={\sf KP}\ome+(V=L)+(\Pi_{1}\mbox{{\rm -Collection}})+(\ome_{1})$ 
{\rm denotes an extension of the Kripke-Platek set theory with the axioms of infinity,
constructibility,
$\Pi_{1}$-Collection
and the following axiom:}
\[
(\ome_{1})\: \:
\exi \kap\fal\alp<\kap\exi\bet,\gam<\kap[\alp<\bet<\gam \land L_{\gam}=rng(F_{\bet\cup\{\kap\}}^{\Sig_{1}})\land 
\mbox{{\rm Hull}}_{\Sig_{1}}(\bet\cup\{\kap\})\cap \kap\subset \bet]
\]
{\rm where}
$F_{\bet\cup\{\kap\}}^{\Sig_{1}}:\mbox{{\rm Hull}}_{\Sig_{1}}(\bet\cup\{\kap\})\to L_{\gam}$
{\rm is the Mostowski collapsing map, and $\mbox{{\rm Hull}}_{\Sig_{1}}(x)$ is the $\Sig_{1}$-Skolem hull
of sets $x$ in the universe.}
\edf

From Theorem \ref{th:cofinalitylocal} we see that
$T_{1}\vdash\exi\kap\fal\alp<\kap(\alp<cf(\kap))$.

\section{A theory equivalent to $T_{1}$}\label{sect:Ztheory}

Referring Theorem \ref{th:cofinalitylocal}  let us interpret $T_{1}$ to another theory.
The base language here is $\{\in\}$.

Let $\rho_{0}$ denotes the least ordinal above the least uncountable ordinal  $\ome_{1}$ such that 
$L_{\rho_{0}}\models (\Pi_{1}\mbox{{\rm -Collection}})$.
$F_{X}(x):=F^{\Sig_{1}}(x;\rho_{0},X)$ and $\mbox{Hull}(X):=\mbox{Hull}_{\Sig_{1}}^{\rho_{0}}(X)$.

The predicate $P$ is intended to denote the relation
$P(x,y)$ iff $x=F_{x\cup\{\ome_{1}\}}(\ome_{1})$ and
$y=F_{x\cup\{\ome_{1}\}}(\rho_{0})$.
Also  the predicate $P_{\rho_{0}}(x)$ is intended to denote the relation
$P_{\rho_{0}}(x)$ iff $x=F_{x}(\rho_{0})$.

\bdf\label{df:regext}
$\mbox{{\rm T}}(\ome_{1})$ {\rm denotes the set theory defined as follows.}
\benu
\item
{\rm Its language is} $\{\in, P,P_{\rho_{0}},\ome_{1}\}$ {\rm for a binary predicate} $P$, {\rm a unary predicate} $P_{\rho_{0}}$
 {\rm  and an individual constant} $\ome_{1}$.

\item
{\rm Its axioms are obtained from those of} $\KP\ome+(\Pi_{1}\mbox{{\rm -Collection}})$ 
{\rm
 in the expanded language
\footnote{
This means that the predicates $P,P_{\rho_{0}}$ do not occur in 
$\Del_{0}$-formulae
for $\Del_{0}$-Separation and 
$\Pi_{1}$-formulae
$\Pi_{1}$-Collection.
},
the axiom of constructibility}
$V=L$
{\rm together with the axiom schema saying that}
$\ome_{1}$ 
 {\rm is an uncountable regular ordinal, cf. (\ref{eq:Z2}) and (\ref{eq:Z1}),
 and if} $P(x,y)$ {\rm then}  $x$ 
 {\rm is a critical point of the} $\Sig_{1}${\rm -elementary embedding from}
$L_{y}\cong \mbox{{\rm Hull}}(x\cup\{\ome_{1}\})$ {\rm to the universe}
$L_{\rho_{0}}${\rm , cf. (\ref{eq:Z1}), and if}
$P_{\rho_{0}}(x)$ {\rm then} $x$ {\rm is a critical point of the} $\Sig_{1}${\rm -elementary embedding from}
$L_{x}\cong \mbox{{\rm Hull}}(x)$ {\rm to the universe}
$L_{\rho_{0}}${\rm , cf.(\ref{eq:Z4}):}
{\rm for a formula} $\vphi$ {\rm and an ordinal} $\alp$,
$\vphi^{\alp}$ {\rm denotes the result of restricting every unbounded quantifier}
$\exi z,\fal z$ {\rm in} $\vphi$ {\rm to} $\exi z\in L_{\alp}, \fal z\in L_{\alp}$.

 \benu
 
 \item
 $x\in Ord$ {\rm is a} $\Del_{0}${\rm -formula saying that `}$x$ {\rm is an ordinal'.}
 \\
$ (\ome<\ome_{1}\in Ord)$, $(P(x,y) \to \{x,y\}\subset Ord \land  x<y<\ome_{1})$
{\rm and}
$(P_{\rho_{0}}(x) \to x\in Ord)$.

 \item
\beqn\label{eq:Z1}
P(x,y) \to a\in L_{x}  \to \vphi[\ome_{1},a] \to \vphi^{y}[x,a]
\eeqn
{\rm for any} $\Sig_{1}${\rm -formula} $\vphi$ {\rm in the language} $\{\in\}$.

.

\item
\beqn\label{eq:Z2}
a\in Ord\cap\ome_{1} \to \exi x, y\in Ord\cap\ome_{1}[a<x\land P(x,y)]
\eeqn

\item
\beqn\label{eq:Z4}
P_{\rho_{0}}(x) \to a\in L_{x} \to \vphi[a] \to \vphi^{x}[a]
\eeqn
{\rm for any} $\Sig_{1}${\rm -formula} $\vphi$ {\rm in the language} $\{\in\}$.

\item
\beqn\label{eq:Z5}
a\in Ord \to \exi x\in Ord[a<x\land P_{\rho_{0}}(x)]
\eeqn

 \eenu
 
\eenu
\edf
{\bf Remark}.
Though the axioms (\ref{eq:Z4}) and (\ref{eq:Z5}) for the $\Pi_{1}$-definable predicate $P_{\rho_{0}}(x)$ are derivable
from $\Pi_{1}$-Collection, the primitive predicate symbol $P_{\rho_{0}}(x)$ is useful for our prof-theoretic study, cf. the proof of Lemma \ref{th:Collapsingthmrho} below.

\blem\label{lem:regularset}
$\mbox{{\rm T}}(\ome_{1})$ is 
a 
conservative extension of
the set theory $T_{1}$.
\elem
\bprf
First consider the axioms of $T_{1}$ in $T(\ome_{1})$.
The axiom $(\ome_{1})$ follows from (\ref{eq:Z1}).
Hence we have shown that $T_{1}$ is contained in $\mbox{{\rm T}}(\ome_{1})$.

Next we show that $\mbox{{\rm T}}(\ome_{1})$ is interpretable in $T_{1}$.
Let $\kap$ be an ordinal in the axiom $(\ome_{1})$.
Interpret the predicate
$P(x,y)\lrarw  \{x,y\}\subset Ord \land 
(\mbox{{\rm Hull}}(x\cup\{\kap\})\cap \kap\subset x)
\land
(y=\sup\{F_{x\cup\{\kap\}}(a):a\in\mbox{{\rm Hull}}(x\cup\{\kap\})\})$.
We see from Theorem \ref{th:cofinalitylocal} that
the interpreted (\ref{eq:Z1}) and (\ref{eq:Z2})
are provable in $T_{1}$.

It remains to show the interpreted (\ref{eq:Z4}) and (\ref{eq:Z5}) in $T_{1}$.
It suffices to show that given an ordinal $\alp$,
there exists an ordinal $x>\alp$
such that $\mbox{Hull}(x)\cap Ord\subset x$.

First we show that for any $\alp$ there exists a $\bet$
such that $\mbox{Hull}(\alp)\cap Ord\subset \bet$.
By Proposition \ref{prp:rhHull}
let $h^{\rho_{0}}_{\Sig_{1}}$ be the $\Del_{1}$-surjection from the $\Sig_{1}$-subset
$dom(h^{\rho_{0}}_{\Sig_{1}})$ of $\ome\times\alp$ to
$\mbox{Hull}(\alp)$, which is 
a $\Sig_{1}$-class.
From $\Sig_{1}$-Separation 
we see that $dom(h^{\rho_{0}}_{\Sig_{1}})$ is a set.
Hence by $\Sig_{1}$-Collection,
$\mbox{Hull}(\alp)=rng(h^{\rho_{0}}_{\Sig_{1}})$ is a set.
Therefore the ordinal $\sup(\mbox{Hull}(\alp)\cap Ord)$ exists
in the universe.

As in Proposition \ref{prp:complexcr}
we see that
$X=\{(\alp,\bet): \bet=\min\{\bet\in Ord: \mbox{{\rm Hull}}(\alp)\cap Ord\subset \bet\}\}$ 
is  a set in $L_{\rho_{0}}$ as follows.
Let $\vphi(\bet)$ be the $\Pi_{1}$-predicate $\vphi(\bet):\Lrarw 
\fal \gam\in Ord[\gam\in\mbox{{\rm Hull}}(\alp)\to \gam\in \bet]$.
Then
$\bet=\min\{\bet: \mbox{{\rm Hull}}(\alp)\cap Ord\subset \bet\}$ iff
$\vphi(\bet)
\land \fal \gam<\bet\lnot \vphi(\gam)$, which is  
$Bool(\Sig_{1}(L_{\rho_{0}}))$ by $\Pi_{0}\mbox{{\rm -Collection}}$.
Hence $X$ is a set in $L_{\rho_{0}}$.

Define recursively ordinals $\{x_{n}\}_{n}$ as follows.
$x_{0}=\alp+1$, and $x_{n+1}$ is defined to be the least ordinal $x_{n+1}$ 
such that
$\mbox{{\rm Hull}}(x_{n})\cap Ord\subset x_{n+1}$, i.e., $(x_{n},x_{n+1})\in X$.
We see inductively that such an ordinal $x_{n}$ exists.
Moreover $n\mapsto x_{n}$ is a $\Del_{1}$-map.
Then $x=\sup_{n}x_{n}<\rho_{0}$ is a desired one.
\eprf

\section{Ordinals for $\ome_{1}$}\label{sect:ordinalnotation}\label{sec:ordinalinacc}
For our proof-theoretic analysis of $T_{1}$, we need to talk about `ordinals' 
less than the next epsilon number
to the order type of the class of ordinals inside $T_{1}$.
Let $Ord^{\varepsilon}\subset V$ and $<^{\varepsilon}$ be $\Delta$-predicates such that
for any transitive and wellfounded model $V$ of $\mbox{{\sf KP}}\omega$,
$<^{\varepsilon}$ is a well ordering of type $\varepsilon_{\rho_{0}+1}$ on $Ord^{\varepsilon}$
for the order type $\rho_{0}$ of the class $Ord$ in $V$.
$<^{\veps}$ is seen to be a canonical ordering as stated in
the following Proposition \ref{prp:canonical}.

\bprp\label{prp:canonical}
\benu
\item\label{prp:canonical0}
$\KP\ome$ proves the fact that $<^{\veps}$ is a linear ordering.

\item\label{prp:canonical2}
For any formula $\vphi$
and each $n<\ome$,
\beqn\label{eq:trindveps}
\KP\ome\vdash\fal x\in Ord^{\veps}(\fal y<^{\veps}x\,\vphi(y)\to\vphi(x)) \to 
\fal x<^{\veps}\ome_{n}(\rho_{0}+1)\vphi(x)
\eeqn
\eenu
\eprp

In what follows of this section we work in $T_{1}$.
For simplicity let us identify the code $x\in Ord^{\veps}$ with
the `ordinal' coded by $x$,
and $<^{\veps}$ is denoted by $<$ when no confusion likely occurs.
Note that the ordinal $\rho_{0}$ 
is the order type of the class of ordinals in the intended model $L_{\rho_{0}}$ of $T_{1}$.
Define simultaneously 
 the classes $\calh_{\alp}(X)\subset L_{\rho_{0}}\cup \veps_{\rho_{0}+1}$
and the ordinals $\Psi_{\ome_{1}} (\alp)$ and $\Psi_{\rho_{0}}(\alp)$ 
for $\alp<^{\veps}\veps_{\rho_{0}+1}$ and sets $X\subset L_{\ome_{1}}$ as follows.
We see that $\calh_{\alp}(X)$ and $\Psi_{\kap} (\alp)\, (\kap\in\{\ome_{1},\rho_{0}\})$ are (first-order) definable as a fixed point in $T_{1}$, cf. Proposition \ref{prp:definability}.

Recall that $\mbox{Hull}(X)=\mbox{Hull}_{\Sig_{1}}^{\rho_{0}}(X)\subset L_{\rho_{0}}$ and
$F_{X}(x)=F^{\Sig_{1}}(x;\rho_{0},X)$ with $F_{X}:\mbox{Hull}(X)\to L_{\gam}$ 
for $X\subset L_{\rho_{0}}$ and a $F_{X}(\rho_{0})=\gam\leq\rho_{0}$.

\bdf\label{df:Cpsiregularsm}

$\calh_{\alp}(X)$ {\rm is the} 
{\rm Skolem hull of} $\{0,\ome_{1},\rho_{0}\}\cup X$ {\rm under the functions} 
$+,
 \alp\mapsto\ome^{\alp}, 
 \Psi_{\ome_{1}}\restrict\alp, \Psi_{\rho_{0}}\restrict \alp${\rm , the}
 $\Sig_{1}${\rm -definability,}
{\rm and the Mostowski collapsing functions}
$(x,d)\mapsto F_{x\cup\{\ome_{1}\}}(d)\,
(\mbox{{\rm Hull}}(x\cup\{\ome_{1}\})\cap\ome_{1}\subset x)$ {\rm and}
$d\mapsto F_{x}(d)\, (\mbox{{\rm Hull}}(x)\cap\rho_{0}\subset x)$.

\benu
\item
$\{0,\ome_{1},\rho_{0}\}\cup X\subset\calh_{\alp}(X)$.

\item
 $x, y \in \calh_{\alp}(X)\Rarw x+ y,\ome^{x}\in \calh_{\alp}(X)$.

\item
$\gam\in \calh_{\alp}(X)\cap\alp
\Rarw 
\Psi_{\kap}(\gam)\in\calh_{\alp}(X)
$
{\rm for} $\kap\in\{\ome_{1},\rho_{0}\}$.

\item
$
\mbox{{\rm Hull}}(\calh_{\alp}(X)\cap L_{\rho_{0}})
\subset \calh_{\alp}(X)
$.

{\rm Namely for any} $\Sig_{1}${\rm -formula} $\vphi[x,\vec{y}]$ {\rm in the language} $\{\in\}$
{\rm and parameters} $\vec{a}\subset \calh_{\alp}(X)\cap L_{\rho_{0}}${\rm , if} 
$b\in L_{\rho_{0}}$,
$(L_{\rho_{0}},\in)\models\vphi[b,\vec{a}]$ {\rm and} 
$(L_{\rho_{0}},\in)\models\exi!x\,\vphi[x,\vec{a}]${\rm , then} $b\in\calh_{\alp}(X)$.

\item
{\rm If} 
$x\in\calh_{\alp}(X)\cap \ome_{1}$ {\rm with} 
$\mbox{{\rm Hull}}(x\cup\{\ome_{1}\})\cap\ome_{1}\subset x$, 
{\rm and}
$d\in (\mbox{{\rm Hull}}(x\cup\{\ome_{1}\})\cap\calh_{\alp}(X))\cup\{\rho_{0}\}$,
{\rm then}
$F_{x\cup\{\ome_{1}\}}(d)\in\calh_{\alp}(X)$.

\item
{\rm If} 
$x\in\calh_{\alp}(X)\cap \rho_{0}$ {\rm with} 
$\mbox{{\rm Hull}}(x)\cap \rho_{0}\subset x$,
{\rm and}
$d\in (\mbox{{\rm Hull}}(x)\cap\calh_{\alp}(X))\cup\{\rho_{0}\}$,
{\rm then}
$F_{x}(d)\in\calh_{\alp}(X)$.

\eenu

{\rm For} $\kap\in\{\ome_{1},\rho_{0}\}$
\[
\Psi_{\kap}(\alp):=
\min\{\bet\leq\kap :  \calh_{\alp}(\bet)\cap \kap \subset \bet\}
.\]

\edf

The ordinal $\Psi_{\kap}(\alp)$ is well defined and $\Psi_{\kap}(\alp)\leq \kap$
for $\kap\in\{\ome_{1},\rho_{0}\}$.

\bprp\label{prp:clshull}
\benu
\item\label{prp:clshull.0}
$\calh_{\alp}(X)$ is closed under $\Sig_{1}$-definability:
$
\vec{a}\subset\calh_{\alp}(X)\cap L_{\rho_{0}}  \Rarw \mbox{{\rm Hull}}(\vec{a})\subset \calh_{\alp}(X)
$.

\item\label{prp:clshull.2}

$
\mbox{{\rm Hull}}(\Psi_{\ome_{1}}(\alp)\cup\{\ome_{1}\})\cap\ome_{1} =\Psi_{\ome_{1}}(\alp)
$
and
$
\mbox{{\rm Hull}}(\Psi_{\rho_{0}}(\alp))\cap\rho_{0} =\Psi_{\rho_{0}}(\alp)>\ome_{1}
$

\item\label{prp:clshull.1}
$\calh_{\alp}(X)$ is closed under the Veblen function $\vphi$ on $\rho_{0}$,
$x, y \in \calh_{\alp}(X)\cap \rho_{0} \Rarw \vphi xy\in \calh_{\alp}(X)$.

\item\label{prp:clshull.3}
If
$x\in\calh_{\alp}(X)\cap\ome_{1}$, $\mbox{{\rm Hull}}(x\cup\{\ome_{1}\})\cap\ome_{1}\subset x$,
 and
$\del\in(\mbox{{\rm Hull}}(x\cup\{\ome_{1}\})\cap \calh_{\alp}(X))\cup\{\rho_{0}\}$,
then
$F_{x\cup\{\ome_{1}\}}(\del)\in\calh_{\alp}(X)$.

\item\label{prp:clshull.4}
If
$x\in\calh_{\alp}(X)\cap\rho_{0}$, $\mbox{{\rm Hull}}(x)\cap\rho_{0}\subset x$,
 and
$\del\in(\mbox{{\rm Hull}}(x)\cap \calh_{\alp}(X))\cup\{\rho_{0}\}$,
then
$F_{x}(\del)\in\calh_{\alp}(X)$.

\eenu
\eprp

The following Proposition \ref{prp:definability} is easy to see.

\bprp\label{prp:definability}
Both of 
$x=\calh_{\alp}(X)$ and
$y=\Psi_{\kap}(\alp)\,(\kap\in \{\ome_{1},\rho_{0}\})$
are $\Sig_{2}$-predicates as fixed points in $\KP\ome$. 
\eprp


\blem\label{lem:lowerbndreg}
For each $n<\ome$,
\[
T_{1}\vdash \fal\alp<\ome_{n+1}(\rho_{0}+1)
\fal\kap\in\{\ome_{1},\rho_{0}\}\exi x<\kap[x=\Psi_{\kap}(\alp)]
.\]
\elem
\bprf
Let $\kap\in\{\ome_{1},\rho_{0}\}$.
By Proposition \ref{prp:definability}
 both 
$x=\calh_{\alp}(\bet)$ and
$y=\Psi_{\kap}(\alp)$ are $\Sig_{2}$-predicates.
We show
that 
$A(\alp) :\Lrarw
 \fal\bet<\rho_{0}\exi x[x=\calh_{\alp}(\bet)] \land \fal\kap\in\{\ome_{1},\rho_{0}\}\exi\bet<\kap[\Psi_{\kap}(\alp)=\bet]$ is 
progressive.
Then $\fal\alp<\ome_{n+1}(\rho_{0}+1)
\fal\kap\in\{\ome_{1},\rho_{0}\}\exi x<\kap[x=\Psi_{\kap}(\alp)]$ will follow from transfinite induction up to $\ome_{n+1}(\rho_{0}+1)$, cf. (\ref{eq:trindveps}) in Proposition \ref{prp:canonical}.

Assume $\fal\gam<\alp\, A(\gam)$ as our IH.
Since $dom(h_{\Sig_{1}}^{\rho_{0}})$ is a $\Sig_{1}$-subset of $\ome\times \bet$ for $\bet<\rho_{0}$, it is a set by $\Sig_{1}$-Separation.
Then so is the image $\mbox{{\rm Hull}}(\bet)$ of the $\Del_{1}$-map $h_{\Sig_{1}}^{\rho_{0}}$.
Hence $\fal \bet<\rho_{0}\exi h[h=\mbox{{\rm Hull}}(\bet)]$.

We see from this, IH and 
$\Sig_{2}$-Collection that $\fal\bet<\rho_{0}\exi x[x=\calh_{\alp}(\bet)=\bigcup_{m}\calh_{\alp}^{m}(\bet)]$,
where $\calh_{\alp}^{m}(\bet)$ is an $m$-th stage of the construction of $\calh_{\alp}(\bet)$ such that
$x=\calh_{\alp}^{m}(\bet)$ is a $\Sig_{2}$-predicate.

Define recursively ordinals $\{\bet_{m}\}_{m}$ for $\kap\in \{\ome_{1},\rho_{0}\}$ as follows.
$\bet_{0}=0$, and 
$\bet_{m+1}$ is defined to be the least ordinal $\bet_{m+1}\leq\kap$ such that
$\calh_{\alp}(\bet_{m})\cap\kap\subset\bet_{m+1}$.

We see inductively that $\bet_{m}<\kap$ using the regularity of $\kap$ and the facts that
$\fal \bet<\kap\exi x[x=\calh_{\alp}(\bet) \land card(x)<\kap]$, where
$card(x)<\kap$ designates that there exists a surjection $f:\gam\to x$ for a $\gam<\kap$ and $f\in L_{\rho_{0}}$.
Moreover $m\mapsto\bet_{m}$ is a $\Sig_{2}$-map.
Therefore $\bet=\sup_{m}\bet_{m}<\kap$ enjoys 
$\calh_{\alp}(\bet)\cap\kap\subset\bet$.
\eprf

\section{Operator controlled derivations for $T_{1}$}\label{sect:controlledOme}

\subsection{An intuitionistic fixed point theory $\mbox{FiX}^{i}(T_{1})$}\label{subsec:intfixZFL}

Let us introduce an intuitionistic fixed point theory $\mbox{FiX}^{i}(T_{1})$ over the set theory $T_{1}$.
Fix an $X$-strictly positive formula $\mathcal{Q}(X,x)$ in the language $\{\in,=,X\}$ with an extra unary predicate symbol $X$.
In $\mathcal{Q}(X,x)$ the predicate symbol $X$ occurs only strictly positive.
The language of $\mbox{FiX}^{i}(T_{1})$ is $\{\in,=,Q\}$ with a fresh unary predicate symbol $Q$.
The axioms in $\mbox{FiX}^{i}(T_{1})$ consist of the following:
\benu
\item
All derivable sentences in $T_{1}$ in the language $\{\in\}$.

\item
Induction schema for any formula $\vphi$ in $\{\in,=,Q\}$:
\\
$
\fal x(\fal y\in x\,\vphi(y)\to\vphi(x))\to\fal x\,\vphi(x)
$.

\item
Fixed point axiom:
$\fal x[Q(x)\lrarw \mathcal{Q}(Q,x)]$.
\eenu

The underlying logic in $\mbox{FiX}^{i}(T_{1})$ is defined to be the intuitionistic first-order predicate logic with equality.

\blem\label{lem:vepsfix}
Let $<^{\veps}$ denote a $\Del_{1}$-predicate defined in section \ref{sec:ordinalinacc}.
For each $n<\ome$ and each formula $\vphi$ in $\{\in,=,Q\}$,
\[
\mbox{{\rm FiX}}^{i}(T_{1})\vdash\fal x(\fal y<^{\veps}x\,\vphi(y)\to\vphi(x)) \to 
\fal x<^{\veps}\ome_{n}(\rho_{0}+1)\vphi(x)
.\]
\elem

The following Theorem \ref{th:consvintfix} is shown in \cite{intfixset}.

\bth\label{th:consvintfix}
$\mbox{{\rm FiX}}^{i}(T_{1})$ is a conservative extension of $T_{1}$.
\end{theorem}

\subsection{Classes of formulae}

In this subsection we work in $T_{1}$.

The language $\calL_{c}$ is obtained from $\{\in, P,P_{\rho_{0}},\ome_{1}\}$ 
by adding names (individual constants) $c_{a}$
of each set $a\in L_{\rho_{0}}$.
$c_{a}$ is identified with $a$.
A {\it term\/} in $\calL_{c}$ is either a variable or a constant in $L_{\rho_{0}}$.
Formulae in this language are defined in the next definition.
Formulae are assumed to be in negation normal form.

\bdf\label{df:OTfml}
\benu
\item
{\rm Let} $t_{1},\ldots,t_{m}$ {\rm be terms.}
{\rm For each} $m${\rm -ary predicate constant} $R\in\{\in,P,P_{\rho_{0}}\}$ 
$R(t_{1},\ldots,t_{m})$ {\rm and} $\lnot R(t_{1},\ldots,t_{m})$ {\rm are formulae,
where} $m=1,2$.
{\rm These are called} literals.


\item
{\rm If} $A$ {\rm and} $B$ {\rm are formulae, then so are} $A\land B$ {\rm and} $A\lor B$.

\item
{\rm Let} $t$ {\rm be a term.}
{\rm If} $A$ {\rm is a formula and the variable} $x$ {\rm does not occur in} $t${\rm , then}
$\exi x\in t\, A$ {\rm and} $\fal x\in t\, A$ {\rm are formulae.}
$\exi x\in t\, A$, $\fal x\in t$ {\rm are} bounded quantifiers.

\item
{\rm If} $A$ {\rm is a formula and} $x$ {\rm a variable, then}
$\exi x\, A$ {\rm and} $\fal x\, A$ {\rm are formulae.}
{\rm Unbounded quantifiers} $\exi x,\fal x$ {\rm are denoted by} $\exi x\in L_{\rho_{0}},\fal x\in L_{\rho_{0}}${\rm , resp.}
\eenu
\edf

For formulae $A$ in $\calL_{c}$, 
$\sfqk(A)$ denotes the finite set 
of sets 
$a\in L_{\rho_{0}}$ 
which are bounds of bounded quantifiers $\exi x\in a,\fal x\in a$
occurring in $A$.
Moreover
$\sfk(A)$ denotes the set of
sets occurring in $A$.
$\sfk(A)$ is defined to include 
bounds of bounded quantifiers.
By definition we set $0\in\sfqk(A)$.
Thus $0\in\sfqk(A)\subset\sfk(A)\subset L_{\rho_{0}}$.


\bdf
\benu
\item
$\sfk(\lnot A)=\sfk(A)$ {\rm and similarly for} $\sfqk$.
\item
$\sfqk(M)=\{0\}$ {\rm for any literal} $M$.
\item
$\sfk(Q(t_{1},\ldots,t_{m}))=(\{t_{1},\ldots,t_{m}\}\cap L_{\rho_{0}})\cup\{0\}$
{\rm for literals} $Q(t_{1},\ldots,t_{m})$ {\rm with predicates} $Q$ {\rm in the set} 
$\{\in,P,P_{\rho_{0}}\}$.
\item
$\sfk(A_{0}\lor A_{1})=\sfk(A_{0})\cup\sfk(A_{1})$ {\rm and similarly for} $\sfqk$.
\item
{\rm For unbounded quantifiers},
$\sfk(\exi x\, A(x))=\sfk(A(x))$ {\rm and similarly for} $\sfqk$.
\item
{\rm For bounded quantifiers with} $a\in L_{\rho_{0}}$,
$\sfk(\exi x\in a\, A(x))=\{a\}\cup\sfk(A(x))$ {\rm and similarly for} $\sfqk$.
\item
{\rm For variables} $y$,
$\sfk(\exi x\in y\, A(x))=\sfk(A(x))$ {\rm and similarly for} $\sfqk$.
\item
{\rm For sets} $\Gam$ {\rm of formulae}
$\sfk(\Gam):=\bigcup\{\sfk(A):A\in\Gam\}$.
\eenu
\edf

For example
$\sfqk(\exi x\in a\, A(x))=\{a\}\cup\sfqk(A(x))$
if $a\in L_{\rho_{0}}$.

\bdf
{\rm For} $a\in L_{\rho_{0}}\cup\{L_{\rho_{0}}\}${\rm ,} $\rk_{L}(a)$ {\rm denotes the} $L$-rank {\rm of} $a$.
\[
\rk_{L}(a):=
\left\{
\begin{array}{ll}
\min\{\alp\in Ord: a\in L_{\alp+1}\} & a\in L_{\rho_{0}}
\\
\rho_{0} & a=L_{\rho_{0}}
\end{array}
\right.
\]
\edf

\bdf\label{df:fmlclasses}
\benu

\item
$A\in\Del_{0}$ {\rm iff there exists a} $\Del_{0}${\rm -formula} $\tht[\vec{x}]$ 
{\rm in the language} $\{\in\}$ {\rm and terms} $\vec{t}$ {\rm in} $\calL_{c}$
 {\rm such that} $A\equiv\tht[\vec{t}]$.
{\rm This means that} $A$ {\rm is bounded, and the predicates} $P,P_{\rho_{0}}$ {\rm do not occur in} $A$.

\item
{\rm Putting} $\Sig_{0}:=\Pi_{0}:=\Del_{0}${\rm , the classes}
$\Sig_{m}$ {\rm and} $\Pi_{m}$ {\rm of formulae in the language} $\calL_{c}$
{\rm are defined as usual, where by definition}
$\Sig_{m}\cup\Pi_{m}\subset\Sig_{m+1}\cap\Pi_{m+1}$.

{\rm Each formula in} $\Sig_{m}\cup\Pi_{m}$ {\rm is in prenex normal form with alternating unbounded quantifiers and}
$\Del_{0}${\rm -matrix.}

\item
{\rm The set} $\Sig^{\Sig_{n}}(\lam)$ {\rm of sentences is defined recursively as follows. Let} $\{a,b,c\}\subset L_{\rho_{0}}$ {\rm and} $d\in L_{\rho_{0}}\cup\{L_{\rho_{0}}\}$.
 \benu
 \item
{\rm Each} $\Sig_{n}${\rm -sentence is in} $\Sig^{\Sig_{n}}(\lam)$.

 \item
 {\rm Each literal including} $Reg(a), P(a,b,c),P_{I,n}(a)$ {\rm and its negation
 is in} $\Sig^{\Sig_{n}}(\lam)$.

 \item
 $\Sig^{\Sig_{n}}(\lam)$ {\rm is closed under propositional connectives} $\lor,\land$.

 \item
{\rm Suppose}
 $\fal x\in d\, A(x)\not\in\Del_{0}${\rm . Then}
 $\fal x\in d\, A(x)\in \Sig^{\Sig_{n}}(\lam)$ {\rm iff} $A(\emptyset)\in \Sig^{\Sig_{n}}(\lam)$ {\rm and}
 $\rk_{L}(d)<\lam$.
 
 \item
{\rm Suppose}
 $\exi x\in d\, A(x)\not\in\Del_{0}${\rm . Then}
 $\exi x\in d\, A(x)\in \Sig^{\Sig_{n}}(\lam)$ {\rm iff} $A(\emptyset)\in \Sig^{\Sig_{n}}(\lam)$ {\rm and}
 $\rk_{L}(d)\leq\lam$.
 \eenu

\item
{\rm For a} $\Sig_{1}${\rm -formula} $A(x)$,
$\exi x\in P_{\rho_{0}}\, A(x)$ {\rm is a}
$\Sig_{1}(P_{\rho_{0}})$-formula.

\eenu
\edf

Note that the predicates $P,P_{\rho_{0}}$ do not occur in $\Sig_{m}$-formulae.

\bdf\label{df:domFfml}
{\rm Let us extend the domain} $dom(F_{x})=\mbox{{\rm Hull}}(x)$
{\rm of the Mostowski collapse to formulae.}
\[
dom(F_{x})=\{A\in\Sig_{1}\cup\Pi_{1}: \sfk(A)\subset\mbox{{\rm Hull}}(x)\}
.\]
{\rm For} $A\in dom(F_{x})$,
$F_{x}" A$ {\rm denotes the result of replacing each constant} 
$c\in L_{\rho_{0}}$ {\rm by} 
$F_{x}(c)${\rm , 
 each unbounded existential quantifier} $\exi z\in L_{\rho_{0}}$ {\rm by} 
 $\exi z\in L_{F_{x}(\rho_{0})}${\rm ,
and each unbounded universal quantifier} $\fal z\in L_{\rho_{0}}$ {\rm by} 
$\fal z\in L_{F_{x}(\rho_{0})}$.

{\rm For sequent, i.e., finite set of sentences} $\Gam\subset dom(F_{x})${\rm , put}
 $F_{x}"\Gam=\{F_{x}" A: A\in\Gam\}$.
 
\edf

The assignment of disjunctions $A\simeq\bigvee(A_{\iota})_{\iota\in J}$ or
conjunctions  $A\simeq\bigwedge(A_{\iota})_{\iota\in J}$
to sentences $A$ is defined as in
\cite{Buchholz} {\it except\/} for $\Sig_{1}\cup\Pi_{1}$-sentences.


\bdf\label{df:assigndc}
\benu
\item\label{df:assigndc0}
{\rm If} $M$ {\rm is one of the literals} $a\in b,a\not\in b${\rm , then for} $J:=0$
\[
M:\simeq
\left\{
\begin{array}{ll}
\bigvee(A_{\iota})_{\iota\in J} & \mbox{{\rm if }} M \mbox{ {\rm is false (in }} $L$\mbox{{\rm )}}
\\
\bigwedge(A_{\iota})_{\iota\in J} &  \mbox{{\rm if }} M \mbox{ {\rm is true}}
\end{array}
\right.
\]

\item
$(A_{0}\lor A_{1}):\simeq\bigvee(A_{\iota})_{\iota\in J}$
{\rm and}
$(A_{0}\land A_{1}):\simeq\bigwedge(A_{\iota})_{\iota\in J}$
{\rm for} $J:=2$.

\item
$
P(b,c) :\simeq
\bigvee(0\not\in 0)_{\iota\in J} 
\mbox{ {\rm and }}
\lnot P(b,c) :\simeq
\bigwedge(0\in 0)_{\iota\in J} 
$
{\rm with}
\[
J:=
\left\{
\begin{array}{ll}
1 & \mbox{{\rm if }}  
 \exi \alp[b=\Psi_{\ome_{1}}(\alp) \spand
c=F_{b\cup\{\ome_{1}\}}(\rho_{0})]
\\
0 & \mbox{{\rm otherwise}}
\end{array}
\right.
.\]

\item
$
P_{\rho_{0}}(a) :\simeq
\bigvee(0\not\in 0)_{\iota\in J} 
\mbox{ {\rm and }}
\lnot P_{\rho_{0}}(a) :\simeq
\bigwedge(0\in 0)_{\iota\in J} 
$
{\rm with}
\[
J:=
\left\{
\begin{array}{ll}
1 & \mbox{{\rm if }}  
 \exi \alp[a=\Psi_{\rho_{0}}(\alp)]
\\
0 & \mbox{{\rm otherwise}}
\end{array}
\right.
.\]

\item
{\rm Let} $\exi z\in b \, \tht[z]\in\Sig_{0}$ {\rm for} $b\in L_{\rho_{0}}\cup\{L_{\rho_{0}}\}$.
{\rm Then for the set}
\[
d:=\mu z\in b \, \tht[z] :=\left\{
\begin{array}{ll}
\min_{<_{L}}\{d : d\in b \land \tht[d]\}  & \mbox{{\rm if }} \exi z\in b \, \tht[z]
\\
0 & \mbox{{\rm otherwise}}
\end{array}
\right.
\]
{\rm with a canonical well ordering} $<_{L}$ {\rm on} $L$, {\rm and}
$J=\{d\}$
\beqnarrs
\exi z\in b\, \tht[z] & :\simeq & \bigvee(d\in b\land \tht[d])_{d\in J}
\\
\fal z\in b\, \lnot\tht[z] & :\simeq & \bigwedge(d\in b \to \lnot\tht[d])_{d\in J}
\eeqnarrs
{\rm where} $d\in b$ {\rm denotes a true literal, e.g.,} $d\not\in d$ {\rm when} $b=L_{\rho_{0}}$.


\item
{\rm For a} $\Sig_{1}(P_{\rho_{0}})${\rm -sentence} $\exi x\in P_{\rho_{0}}\, A(x)$,
\beqnarrs
\exi x\in P_{\rho_{0}}\, A(x) & \simeq & \bigvee(A(a))_{a\in J} 
\\
\fal x\in P_{\rho_{0}}\, \lnot A(x) & \simeq & \bigwedge(\lnot A(a))_{a\in J} 
\\
\mbox{ {\rm with  }} J & = & \{a: \exi\alp(a=\Psi_{\rho_{0}}(\alp))\}
\eeqnarrs

\item
{\rm Otherwise set for} $a\in L_{\rho_{0}}\cup\{L_{\rho_{0}}\}$ {\rm and} $J:=\{b: b\in a\}$
\[
\exi x\in a\, A(x):\simeq\bigvee(A(b))_{b\in J}
\mbox{ {\rm and }}
\fal x\in a\, A(x):\simeq\bigwedge(A(b))_{b\in J}
.\]

\eenu
\edf

The rank $\rk(A)$ of sentences $A$ is defined by recursion on the  number of symbols occurring in 
$A$.

\bdf\label{df:rank}
\benu
\item\label{df:rank1}
$\rk(\lnot A):=\rk(A)$.

\item\label{df:rank2}
$\rk(a\in b):=0$.

\item\label{df:rank5-1}
$\rk(P(b,c)):=\rk(P_{\rho_{0}}(a)):=1$.


\item\label{df:rank6}
$\rk(A_{0}\lor A_{1}):=\max\{\rk(A_{0}),\rk(A_{1})\}+1$.


\item\label{df:rank7}
$
\rk(\exi x\in a\, A(x)):=
\max\{\ome\alp, \rk(A(\emptyset))+1\}
$
{\rm for} $\alp=\rk_{L}(a)$.

\item\label{df:rank8}
$\rk(\exi x\in P_{\rho_{0}}\, A(x))=\rho_{0}$.
\eenu
\edf

\bprp\label{lem:rank}
Let $A\simeq\bigvee(A_{\iota})_{\iota\in J}$ or $A\simeq\bigwedge(A_{\iota})_{\iota\in J}$.
\benu
\item\label{lem:rank15-1}
$\fal\iota\in J(\sfk(A_{\iota})\subset\sfk(A)\cup\{\iota\})$.

\item\label{lem:rank15}
$A\in \Sig^{\Sig_{n}}(\lam)\Rarw\fal\iota\in J(A_{\iota}\in \Sig^{\Sig_{n}}(\lam))$.

\item\label{prp:rksig3}
For an ordinal $\lam\leq \rho_{0}$ with $\ome\lam=\lam$,
$
\rk(A)<\lam \Rarw A\in\Sig^{\Sig_{n}}(\lam)
$.


\item\label{lem:rank0}
$\rk(A)<\rho_{0}+\ome$.

\item\label{lem:rank1}
$\rk(A)\in\{\ome \, \rk_{L}(a)+i : 
a\in\sfqk(A)\cup\{\rho_{0}\}, 
i\in\ome\}
\subset\mbox{{\rm Hull}}(\sfk(A))$.

\item\label{lem:rank2}
$\fal\iota\in J(\rk(A_{\iota})<\rk(A))$.

\eenu
\eprp

\subsection{Operator controlled derivations}\label{subsec:opcontderivation}

In the remaining parts of this section we work in the intuitionistic fixed point theory $\mbox{FiX}^{i}(T_{1})$.

Sequents are finite sets of sentences, and inference rules are formulated in one-sided sequent calculus.
In what follows by an operator we mean an $\calh_{\gam}[\Tht]$ for a finite set $\Tht$ of sets.

\bdf\label{df:controlderreg}
{\rm Define a relation} $(\calh,\kap)\vdash^{a}_{b}\Gam$ {\rm as follows.}

$(\calh,\kap)\vdash^{a}_{b}\Gam$ {\rm holds if} 
\beqn
\label{eq:controlder}
\{a\}\cup\sfk(\Gam)\subset\calh:=\calh(\emptyset)
\eeqn
{\rm and one of the following
cases holds:}


\bdes
\item[$(\bigvee)$]
$A\simeq\bigvee\{A_{\iota}: \iota\in J\}$, $A\in\Gam$ {\rm and there exist} $\iota\in J$
{\rm and}
 $a(\iota)<a$ {\rm such that}
\beqn\label{eq:bigveebnd}
\rk_{L}(\iota)<\kap \Rarw \rk_{L}(\iota)< a
\eeqn
{\rm and}
$(\calh,\kap)\vdash^{a(\iota)}_{b}\Gam,A_{\iota}$.

\item[$(\bigwedge)$]
$A\simeq\bigwedge\{A_{\iota}: \iota\in J\}$, $A\in\Gam$ {\rm and for every}
$\iota\in J$ {\rm there exists an} $a(\iota)<a$ 
{\rm such that}
$(\calh[\{\iota\}],\kap)\vdash^{a(\iota)}_{b}\Gam,A_{\iota}$.

\item[$(cut)$]
{\rm There exist} $a_{0}<a$ {\rm and} 
$C$
{\rm such that} $\mbox{{\rm rk}}(C)<b$
{\rm and}
$(\calh,\kap)\vdash^{a_{0}}_{b}\Gam,\lnot C$
{\rm and}
$(\calh,\kap)\vdash^{a_{0}}_{b}C,\Gam$.

\item[$(\mbox{P})$]
{\rm There exists} $\alp<\ome_{1}$ {\rm such that} 
$(\exi x,y<\ome_{1}[\alp<x \land P(x,y)])\in\Gam$.

\item[$(\mbox{F}_{x\cup\{\ome_{1}\}})$] 

$x=\Psi_{\ome_{1}}(\bet)\in\calh$ {\rm for a} $\bet$
{\rm and there exist} $a_{0}<a$, 
$\Gam_{0}\subset\Sig_{1}$ {\rm and} $\Lam$ {\rm such that}
$\sfk(\Gam_{0})\subset\mbox{{\rm Hull}}(x\cup\{\ome_{1}\})$,
$\Gam=\Lam\cup (F_{x\cup\{\ome_{1}\}}"\Gam_{0})$
{\rm and}
$
(\calh,\kap)\vdash^{a_{0}}_{b}\Lam,\Gam_{0}$,
{\rm where} $F_{x\cup\{\ome_{1}\}}$ {\rm denotes the Mostowski collapse}
$F_{x\cup\{\ome_{1}\}}: \mbox{{\rm Hull}}(x\cup\{\ome_{1}\})\lrarw L_{F_{x\cup\{\ome_{1}\}}(\rho_{0})}$.

\item[$(\mbox{P}_{\rho_{0}})$]
{\rm There exists} $\alp<\rho_{0}$ {\rm such that} 
$(\exi x<\rho_{0}[\alp<x \land P_{\rho_{0}}(x)])\in\Gam$.

\item[$(\mbox{F}_{x})$] 
$
x=\Psi_{\rho_{0}}(\bet)\in\calh
$ {\rm for a} $\bet$ {\rm and there exist} $a_{0}<a$, $\Gam_{0}\subset\Sig_{1}$
{\rm and} $\Lam$ {\rm such that}
$\sfk(\Gam_{0})\subset\mbox{{\rm Hull}}(x)$,
$\Gam=\Lam\cup (F_{x}"\Gam_{0})$ {\rm and}
$(\calh,\kap)\vdash^{a_{0}}_{b}\Lam,\Gam_{0}$,
{\rm where} $F_{x}$ {\rm denotes the Mostowski collapse}
$F_{x}: \mbox{{\rm Hull}}(x)\lrarw L_{F_{x}(\rho_{0})}$.


\item[$(Ref)$]
$b\geq\rho_{0}$,
{\rm and there exist an ordinal} $a_{0}<a$, {\rm a set} $c$ {\rm and a} $\Sig_{1}(P_{\rho_{0}})${\rm -formula}
$A(x)$ {\rm such that} 
$(\calh,\kap)\vdash^{a_{0}}_{b}\Gam, \fal x\in c\, A(x)$
{\rm and}
$(\calh,\kap)\vdash^{a_{0}}_{b} \fal y\exi x\in c\, \lnot A^{(y)}(x),\Gam$, {\rm where for}
$A(x)\equiv(\exi z\in P_{\rho_{0}}\exi w\, B(x))\, (B\in\Del_{0})$,
$A^{(y)}(x):\equiv(\exi z\in P_{\rho_{0}}\cap y\exi w\in y\, B)$.
\edes
\edf


\blem\label{lem:tautology}{\rm (Tautology)}
If $\sfk(\Gam\cup\{A\})\subset\calh$
then
$(\calh,\rho_{0})\vdash^{2\footnotesize{\rk}(A)}_{0}\Gam,\lnot A, A$.
\elem



\blem\label{lem:embedfund}
Let  $\rk(\fal x\in b\,\vphi[x,c])\leq \rho_{0}+m$ for an $m\geq 1$, 
and $\Tht_{c}=\{\lnot\fal y(\fal x\in y\,\vphi[x,c]\to\vphi[y,c])\}$.
Then for any operator $\calh$, and any $a,c$,
\\
$(\calh[\{c,a\}],\rho_{0})\vdash^{\rho_{0}+2m+2+2\footnotesize{\rk}_{L}(a)}_{\rho_{0}+m+1}\Tht_{c},\fal x\in a\,\vphi[x,c]$.
\elem
Let
\beqnarrs
(\calh,\rho_{0})\vdash_{c}^{<\alp}\Gam & :\Lrarw & \exi\bet<\alp[(\calh,\rho_{0})\vdash_{c}^{\bet}\Gam]
\\
(\calh,\rho_{0})\vdash_{<c}^{<\alp}\Gam & :\Lrarw & \exi d<c[(\calh,\rho_{0})\vdash_{d}^{<\alp}\Gam]
\eeqnarrs

\blem\label{th:embedreg}
Let $A$ be an axiom in $\mbox{{\rm T}}(\ome_{1})$ except Foundation axiom schema and $\Pi_{1}$-Collection.
Then
 $(\calh,\rho_{0})\vdash_{0}^{<\rho_{0}+\ome}A$ for any operator $\calh=\calh_{\gam}$.
\elem

\blem\label{th:embedregthm}{\rm (Embedding)}\\
If $\mbox{{\rm T}}(\ome_{1})\vdash \Gam[\vec{x}]$, there are $m,k<\ome$ such that
for any $\vec{a}\subset L_{\rho_{0}}$,
 $(\calh[\vec{a}],\rho_{0})\vdash_{\rho_{0}+m}^{\rho_{0}\cdot 2+k}\Gam[\vec{a}]$ for any operator 
 $\calh=\calh_{\gam}$.
\elem
\bprf

By Lemma \ref{lem:embedfund}
we have
$(\calh,\rho_{0})\vdash^{\rho_{0}\cdot 2}_{\rho_{0}+m+1}\fal u,z(\fal y(\fal x\in y\,\vphi[x,z]\to\vphi[y,z])\to \vphi[u,z])$
for some $m$.
By Lemmata \ref{lem:tautology} and \ref{th:embedreg} it remains to consider instances 
\[
\fal u\in a\exi v\fal w\, \tht \to \exi c\fal u\in a\exi v\in c\fal w\, \tht
\]
 of
$\Pi_{1}$-Collection, where $\tht\equiv\tht(u,v,w)$ is a $\Del_{0}$-formula in the language $\{\in\}$.

First by Lemma \ref{th:embedreg} with axioms (\ref{eq:Z4}) and (\ref{eq:Z5}) we have
\[
(\calh,\rho_{0})\vdash_{\rho_{0}+1}^{\rho_{0}+\ome}\fal w\,\tht(u,v,w)\lrarw \exi x\in P_{\rho_{0}} \tau(x,u,v)
\]
where $\tau(x,u,v)\equiv[u,v\in L_{x} \land \fal w\in L_{x}\, \tht(u,v,w)]$.
Hence 
\[
(\calh,\rho_{0})\vdash_{<\rho_{0}+\ome}^{<\rho_{0}+\ome\cdot 2}\lnot\fal u\in a\exi v\fal w\, \tht, \fal u\in a\exi x\in P_{\rho_{0}}\exi v\, \tau(x,u,v)
\]
On the other hand we have by Lemma \ref{lem:tautology}
\[
(\calh,\rho_{0})\vdash_{0}^{<\rho_{0}+\ome}\lnot \exi c\fal u\in a\exi x\in P_{\rho_{0}}\cap c\exi v\in c\,\tau,
\exi c\fal u\in a\exi x\in P_{\rho_{0}}\cap c\exi v\in c\,\tau
\]
Hence by the inference $(Ref)$ for the $\Sig_{1}(P_{\rho_{0}})$-formula $\exi x\in P_{\rho_{0}}\exi v\, \tau(x,u,v)$, we obtain
\[
(\calh,\rho_{0})\vdash_{<\rho_{0}+\ome}^{<\rho_{0}+\ome\cdot 2}\lnot\fal u\in a\exi v\fal w\, \tht,\exi c\fal u\in a\exi x\in P_{\rho_{0}}\cap c\exi v\in c\,\tau
\]
Therefore
$(\calh,\rho_{0})\vdash_{<\rho_{0}+\ome}^{\rho_{0}+\ome\cdot 2}\fal u\in a\exi v\fal w\, \tht \to \exi c\fal u\in a\exi v\in c\fal w\, \tht$.
\eprf
\\

In the following Lemma \ref{lem:predcereg}, note that 
$\rk(\exi x<\ome_{1}\exi y<\ome_{1}[\alp<x \land P(x,y)])=\ome_{1}+1$, and
$\rk(\exi x<\rho_{0}[\alp<x \land P_{\rho_{0}}(x)])=\rho_{0}$.

\blem\label{lem:predcereg}{\rm (Predicative Cut-elimination)}
\benu
\item\label{lem:predcereg2}
If 
$(\calh,\kap)\vdash^{b}_{c+\ome^{a}}\Gam
\spand 
[c,c+\ome^{a}[\cap \{\ome_{1}+1,\rho_{0}\}=\emptyset
\spand a\in\calh
\Rarw (\calh,\kap)\vdash^{\vphi ab}_{c}\Gam$.

\item\label{lem:predcereg4}
If 
$(\calh_{\gam},\kap)\vdash^{b}_{\ome_{1}+2}\Gam \spand \gam\in\calh_{\gam}
 \Rarw 
(\calh_{\gam+b},\kap)\vdash^{\ome^{b}}_{\ome_{1}+1}\Gam$.

\item\label{lem:predcereg5}
If 
$(\calh_{\gam},\kap)\vdash^{b}_{\rho_{0}+1}\Gam \spand \gam\in\calh_{\gam}
 \Rarw 
(\calh_{\gam+b},\kap)\vdash^{\ome^{b}}_{\rho_{0}}\Gam$.

\eenu
\elem

{\rm For a formula} $\exi x\in d\, A(x)$ {\rm and ordinals} $\lam=\rk_{L}(d), \alp$,
$(\exi x\in d\, A)^{(\exi\lam\restrict\alp)}$ {\rm denotes the result of restricting the} outermost existential quantifier 
$\exi x\in d$ {\rm to} $\exi x\in L_{\alp}$,
$(\exi x\in d\, A)^{(\exi\lam\restrict\alp)}\equiv
(\exi x\in L_{\alp}\, A)$.

\blem\label{lem:boundednessreg}{\rm (Boundedness)}
Let $\lam\in \{\ome_{1},\rho_{0}\}$, $C\equiv(\exi x\in d\, A)\in\Sig^{\Sig_{2}}(\lam)$ and $C\not\in\{\exi x<\ome_{1} \exi y<\ome_{1}[\alp<x \land P(x,y)]: \alp<\ome_{1}\}\cup\{\exi x<\rho_{0}[\alp<x \land P_{\rho_{0}}(x)]:\alp<\rho_{0}\}$.
\benu
\item\label{lem:boundednessregexi}
$(\calh,\lam)\vdash^{a}_{c}\Lam, C \spand a\leq b\in\calh\cap\lam
\Rarw (\calh,\lam)\vdash^{a}_{c}\Lam,C^{(\exi\lam\restrict b)}$.

\item\label{lem:boundednessregfal}
$(\calh,\kap)\vdash^{a}_{c}\Lam,\lnot C \spand b\in\calh\cap\lam 
\Rarw (\calh,\kap)\vdash^{a}_{c}\Lam,\lnot (C^{(\exi\lam\restrict b)})$.
\eenu
\elem

\blem\label{lem:boundednessrho}{\rm (Boundedness for $\Sig_{1}(P_{\rho_{0}})$)}
\\
Let $C$ be a $\Sig_{1}(P_{\rho_{0}})$-sentence. Then
$
(\calh,\rho_{0})\vdash^{a}_{c}\Lam, C \spand a\leq b\in\calh\cap\rho_{0}
\Rarw (\calh,\rho_{0})\vdash^{a}_{c}\Lam, C^{(L_{b})}
$.
\elem
\bprf
$C^{L_{b}}\equiv (\exi z\in P_{\rho_{0}}\cap L_{b}\exi w\in L_{b}\, B)$ when $C\equiv(\exi z\in P_{\rho_{0}}\exi w\, B)$ with a $\Del_{0}$-formula $B$.
The lemma is seen from (\ref{eq:bigveebnd}).
\eprf

\subsection{Collapsing derivations}\label{subsec:collapse}
In this subsection derivations of $\Sig^{\Sig_{2}}(\ome_{1})$ sentences are shown to be collapsed to derivations
with heights and cut ranks$<\ome_{1}$.

\blem\label{th:Collapsingthmomega1}{\rm (Collapsing below $\ome_{1}$)}\\

Suppose $\gam\in\calh_{\gam}[\Tht]$ with $\Tht\subset\calh_{\gam}(\Psi_{\ome_{1}}(\gam))$, and
$
\Gam\subset\Sig^{\Sig_{2}}(\ome_{1})
$.

Then for $b=\Psi_{\ome_{1}}(\gam+\ome^{\ome_{1}+a})$,
\[
 (\calh_{\gam}[\Tht],\ome_{1})\vdash^{a}_{\ome_{1}+1}\Gam \Rarw
  (\calh_{\gam+\ome^{\ome_{1}+a}+1}[\Tht],\ome_{1})\vdash^{b}_{b} \Gam.
\]
\elem

\blem\label{th:Collapsingthmrho}{\rm (Collapsing below $\rho_{0}$)}\\

Suppose $\gam\in\calh_{\gam}[\Tht]$ with $\Tht\subset\calh_{\gam}(\Psi_{\rho_{0}}(\gam))$, and 
$\Gam\subset\Sig^{\Sig_{2}}(\rho_{0})\cup\Sig_{1}(P_{\rho_{0}})$.

Then for $\hat{a}=\gam+\ome^{\rho_{0}+a}$ 
\[
 (\calh_{\gam}[\Tht],\rho_{0})\vdash^{a}_{\rho_{0}}\Gam \Rarw
  (\calh_{\hat{a}+1}[\Tht],\rho_{0})\vdash^{\Psi_{\rho_{0}}(\hat{a})}_{\Psi_{\rho_{0}}(\hat{a})}\Gam.
\]

\elem
{\bf Proof} by induction on $a$, cf. Lemma \ref{lem:vepsfix}.

First note that $\Psi_{\rho_{0}}(\hat{a})\in\calh_{\hat{a}+1}[\Tht]$
since 
$\hat{a}=\gam+\ome^{\rho_{0}+a}\in\calh_{\gam}[\Tht]\subset\calh_{\hat{a}+1}[\Tht]$
 by the assumption,
$\{\gam,a\}\subset\calh_{\gam}[\Tht]$.
 
Assume $(\calh_{\gam}[\Tht][\Tht_{0}],\rho_{0})\vdash^{a_{0}}_{\rho_{0}}\Gam_{0}$ with 
$\Tht_{0}\subset\calh_{\gam}(\Psi_{\rho_{0}}(\gam))$.
Then by $\gam\leq\hat{a}$, we have 
$\widehat{a_{0}}\in \calh_{\gam}[\Tht][\Tht_{0}]\subset \calh_{\gam}(\Psi_{\rho_{0}}(\gam))\subset\calh_{\hat{a}}(\Psi_{\rho_{0}}(\hat{a}))$.
This yields
that
\[
a_{0}<a \Rarw \Psi_{\rho_{0}}(\widehat{a_{0}})<\Psi_{\rho_{0}}(\hat{a})
\]

Second observe that $\sfk(\Gam)\subset\calh_{\gam}[\Tht]\subset\calh_{\hat{a}+1}[\Tht]$ by $\gam\leq\hat{a}+1$.
 
 Third we have
\[
\sfk(\Gam)\subset\calh_{\gam}(\Psi_{\rho_{0}}(\gam))
\]
When $\Gam$ is one of axioms $(\mbox{{\bf P}})$ and $(\mbox{{\bf P}}_{\rho_{0}})$, there is nothing to show.

Consider the case when the last inference is a $(Ref)$. 
\[
\infer[(Ref)]{(\calh_{\gam}[\Tht],\rho_{0})\vdash^{a}_{\rho_{0}}\Gam}
{
(\calh_{\gam}[\Tht],\rho_{0})\vdash^{a_{0}}_{\rho_{0}}\Gam,\fal x\in c\, A(x)
&
(\calh_{\gam}[\Tht],\rho_{0})\vdash^{a_{0}}_{\rho_{0}}\fal y\exi x\in c\, \lnot A^{(y)}(x),\Gam
}
\]
where $a_{0}<a$ and
$A(x)\equiv(\exi z\in P_{\rho_{0}}\exi w\, B(x))$ is a $\Sig_{1}(P_{\rho_{0}})$-formula with a $\Del_{0}$-formula $B$.

For each $d\in c$ we have by Inversion
\[
(\calh_{\gam}[\Tht\cup\{d\}],\rho_{0})\vdash^{a_{0}}_{\rho_{0}}\Gam, A(d)
\]
where $c\in\calh_{\gam}(\Psi_{\rho_{0}}(\gam))$.
Hence $\rk_{L}(d)<\rk_{L}(c)\in\calh_{\gam}(\Psi_{\rho_{0}}(\gam))\cap\rho_{0}\subset\Psi_{\rho_{0}}(\gam)$,
and $\rk_{L}(d)<\Psi_{\rho_{0}}(\gam)$.
Therefore $d\in \calh_{\gam}(\Psi_{\rho_{0}}(\gam))$.
By IH we have for $\widehat{a_{0}}=\gam+\ome^{\rho_{0}+a_{0}}$ and 
$\bet_{0}=\Psi_{\rho_{0}}(\widehat{a_{0}})\in\calh_{\widehat{a_{0}}+1}[\Tht]$
\[
(\calh_{\widehat{a_{0}}+1}[\Tht\cup\{d\}],\rho_{0})\vdash^{\bet_{0}}_{\bet_{0}}
\Gam,A(d)
\]
Boundedness lemma \ref{lem:boundednessrho} yields
\[
(\calh_{\widehat{a_{0}}+1}[\Tht\cup\{d\}],\rho_{0})\vdash^{\bet_{0}}_{\bet_{0}}
\Gam,A^{(L_{\bet_{0}})}(d)
\]
Since $d\in c$ is arbitrary, we obtain by $(\bigwedge)$
\beqn\label{eq:Collapsingthmrho1}
(\calh_{\widehat{a_{0}}+1}[\Tht],\rho_{0})\vdash^{\bet_{0}+1}_{\bet_{0}}
\Gam,\fal x\in c\, A^{(L_{\bet_{0}})}(x)
\eeqn
On the other hand we have by Inversion for $L_{\bet_{0}}\in\calh_{\widehat{a_{0}}+1}[\Tht]$
\[
(\calh_{\widehat{a_{0}}+1}[\Tht],\rho_{0})\vdash^{a_{0}}_{\rho_{0}}\exi x\in c\, \lnot A^{(L_{\bet_{0}})}(x),\Gam
\]
Since $\exi x\in c\, \lnot A^{(L_{\bet_{0}})}(x)\in\Sig^{\Sig_{2}}(\rho_{0})$, IH yields for
$\widehat{a_{1}}=\widehat{a_{0}}+1+\ome^{\rho_{0}+a_{0}}=\gam+\ome^{\rho_{0}+a_{0}}\cdot 2$
and $\bet_{1}=\Psi_{\rho_{0}}(\widehat{a_{1}})$
\beqn\label{eq:Collapsingthmrho2}
(\calh_{\widehat{a_{1}}+1}[\Tht],\rho_{0})\vdash^{\bet_{1}}_{\bet_{1}}\exi x\in c\, \lnot A^{(L_{\bet_{0}})}(x),\Gam
\eeqn
We have $\rk(\fal x\in c\, A^{(L_{\bet_{0}})}(x))\in\mbox{Hull}(\sfk(\fal x\in c\, A^{(L_{\bet_{0}})}(x)))\cap\rho_{0}\subset\calh_{\widehat{a_{0}}+1}[\Tht]\cap\rho_{0}\subset\calh_{\widehat{a_{0}}+1}(\Psi_{\rho_{0}}(\gam))\cap\rho_{0}\subset\Psi_{\rho_{0}}(\hat{a})$
by Proposition \ref{lem:rank}.\ref{lem:rank1}.

By a $(cut)$ with (\ref{eq:Collapsingthmrho1}) and (\ref{eq:Collapsingthmrho2}) we obtain with 
$\Psi_{\rho_{0}}(\hat{a})>\bet_{1}>\bet_{0}$
\[
(\calh_{\hat{a}+1}[\Tht],\rho_{0})\vdash^{\Psi_{\rho_{0}}(\hat{a})}_{\Psi_{\rho_{0}}(\hat{a})}\Gam
\]

Other case ae seen as in \cite{liftupZF}.
\eprf

\section{{\bf Proof} of Theorem \ref{th:mainthZ}}\label{sect:proof}

For a sentence $\exi x\in L_{\ome_{1}}\vphi$ with a $\Sig_{2}$-formula $\vphi$ in the language $\{\in,\ome_{1}\}$,
assume
$T_{1}\vdash\exi x\in L_{\ome_{1}}\,\vphi$.
Then by Lemmata \ref{lem:regularset} and \ref{th:embedregthm}, pick an $m>0$
such that the fact
$(\calh_{0},\rho_{0})\vdash^{\rho_{0}\cdot 2+m}_{\rho_{0}+m}\exi x\in L_{\ome_{1}}\vphi$ is provable in 
$\mbox{FiX}^{i}(T_{1})$.
In what follows work in $\mbox{FiX}^{i}(T_{1})$.
Predicative Cut Elimination \ref{lem:predcereg}.\ref{lem:predcereg2} and \ref{lem:predcereg}.\ref{lem:predcereg5}
yields
\[
(\calh_{\gam},\rho_{0})\vdash^{a}_{\rho_{0}}\exi x\in L_{\ome_{1}}\vphi
\]
for
$\gam=\ome_{m-1}(\rho_{0}\cdot 2+m)$ and $a=\ome_{m}(\rho_{0}\cdot 2+m)$.
Then Collapsing below $\rho_{0}$ \ref{th:Collapsingthmrho} yields
\[
(\calh_{\ome_{m+1}(\rho_{0}\cdot 2+m)+1},\rho_{0})\vdash^{\bet}_{\bet}\exi x\in L_{\ome_{1}}\vphi
\]
for $\gam+\ome^{\rho_{0}+a}=\ome_{m+1}(\rho_{0}\cdot 2+m)$ and 
$\bet=\Psi_{\rho_{0}}(\ome_{m+1}(\rho_{0}\cdot 2+m))$.
Predicative Cut Elimination \ref{lem:predcereg}.\ref{lem:predcereg2} and \ref{lem:predcereg}.\ref{lem:predcereg4} yields
\[
(\calh_{\ome_{m+1}(\rho_{2}\cdot 2+m)+\vphi\bet\bet},\rho_{0})\vdash^{\vphi\bet\bet}_{\ome_{1}+1}\exi x\in L_{\ome_{1}}\vphi
\]
for 
$\ome_{1}+2+\ome^{\bet}=\bet$ and $\ome^{\vphi\bet\bet}=\vphi\bet\bet$.
A fortiori,
\[
(\calh_{\ome_{m+1}(\rho_{2}\cdot 2+m)+\vphi\bet\bet},\ome_{1})\vdash^{\vphi\bet\bet}_{\ome_{1}+1}\exi x\in L_{\ome_{1}}\vphi
\]
Then Collapsing below $\ome_{1}$ \ref{th:Collapsingthmomega1} yields
\[
(\calh_{\ome_{m+1}(\rho_{2}\cdot 2+m)+(\vphi\bet\bet\cdot) 2+1},\ome_{1})\vdash^{\del}_{\del}\exi x\in L_{\ome_{1}}\vphi
\]
for
$\ome_{m+1}(\rho_{2}\cdot 2+m)+\vphi\bet\bet+\ome^{\ome_{1}+\vphi\bet\bet}+1=\ome_{m+1}(\rho_{2}\cdot 2+m)+(\vphi\bet\bet\cdot) 2+1$
and $\del=\Psi_{\ome_{1}}(\ome_{m+1}(\rho_{2}\cdot 2+m)+(\vphi\bet\bet\cdot) 2)$.

Boundedness \ref{lem:boundednessreg}.\ref{lem:boundednessregexi} yields for $\del<\Psi_{\ome_{1}}(\ome_{n}(\rho_{0}+1))$ with
$n=m+2$
\[
(\calh_{\ome_{n}(\rho_{0}+1)+1},\ome_{1})\vdash^{\del}_{\del}
\exi x\in L_{\Psi_{\ome_{1}}(\ome_{n}(\rho_{0}+1))}\vphi
\]
We see then by transfinite induction up to the countable ordinal $\del$ that
inference rules in the controlled derivation of $\exi x\in L_{\Psi_{\ome_{1}}(\ome_{n}(\rho_{0}+1))}\vphi$ with cut rank$<\ome_{1}$
are $(\bigvee)$, $(\bigwedge)$, $(cut)$, 
and $({\bf F}_{x\cup\{\ome_{1}\}})$,
and since these inference rules are truth-preserving, we conclude again 
by transfinite induction up to $\del$ that
$\exi x\in L_{\Psi_{\ome_{1}}(\ome_{n}(\rho_{0}+1))}\vphi$ is true.

Since the whole proof is formalizable in $\mbox{FiX}^{i}(T_{1})$, 
we conclude
$\mbox{FiX}^{i}(T_{1})\vdash  \exi x\in L_{\Psi_{\ome_{1}}(\ome_{n}(\rho_{0}+1))} \vphi$.
Finally Theorem \ref{th:consvintfix} yields
$T_{1}\vdash  \exi x\in L_{\Psi_{\ome_{1}}(\ome_{n}(\rho_{0}+1))} \vphi$.
This completes a proof of Theorem \ref{th:mainthZ}.

\end{document}